\documentclass[12pt]{article}
\usepackage{amssymb,latexsym,start2e,graphs2,pstcol,pst-plot,amsthm,amsmath}

\newtheorem{thm}{Theorem}[section]

\newtheorem{question}[thm]{Question}

\begin{document}
\pagestyle{plain}

\title{Compositions inside a rectangle and unimodality
}
\author{
Bruce E. Sagan\\[-5pt]
\small Department of Mathematics, Michigan State University,\\[-5pt]
\small East Lansing, MI 48824-1027, USA, \texttt{sagan@math.msu.edu}
}

\date{\today\\[10pt]
	\begin{flushleft}
	                                       \\[5pt]
	\end{flushleft}}

\maketitle

\begin{abstract}
Let $c^{k,l}(n)$ be the number of compositions (ordered
partitions) of the integer $n$ whose Ferrers diagram fits inside a
$k\times l$ rectangle.  The purpose of this note is to give  a simple, algebraic proof of a conjecture
of Vatter that the sequence $c^{k,l}(0),c^{k,l}(1),\ldots,c^{k,l}(kl)$ is unimodal.  The
problem of giving a combinatorial proof of this fact is discussed, but
is still open.
\end{abstract}

\section{Introduction}

Let $\bbN$ and $\bbP$ denote the nonnegative and positive integers,
respectively.  A {\it partition\/} of $n\in\bbN$ is a weakly
decreasing sequence $\la=(\la_1,\ldots,\la_r)$ of positive integers
called {\it parts\/} such that $\sum_i \la_i=n$.  We write $\la\ptn n$
or $|\la|=n$ if $\la$ partitions $n$.  We will also use the  notation
$\la=(n^{m_n},\ldots,1^{m_1})$ where $m_i$ is the number of times $i$
appears as a part in $\la$.  
If $m_i=1$ then the exponent is suppressed and if $m_i=0$ then so is
the base.
The {\it Ferrers diagram\/} of $\la$,
also denoted $\la$,
consists of left-justified rows of squares with $\la_i$ squares in row
$i$.  The Ferrers diagram of $\la=(4,3,3,1)=(4,3^2,1)$ is shown in
Figure~\ref{laka}.

Partitions can be  ordered by letting $\la\le\mu$ if the
Ferrers diagram for $\la$ is contained in the upper left corner of the
one for $\mu$.  Equivalently, $\la_i\le\mu_i$ for all $i$ where we set
$\la_i=0$ if $i$ is greater than the number of parts of $\la$ and
similarly for $\mu$.  The set of partitions under this partial order
is called {\it Young's lattice\/}.  More information about partitions
and this lattice can be found in the books of Andrews~\cite{and:tp}, Sagan~\cite{sag:sym}, or
Stanley~\cite{sta:ec2}.

We say that $\la$ {\it fits inside a $k\times l$ rectangle\/} if
$\la\le(l^k)$.  In other words, $\la$ has at most $k$ parts each of
size at most $l$.  Let $p^{k,l}(n)$ denote the number of such
$\la$ where $\la\ptn n$.  A sequence $a_0,a_1,\ldots,a_r$ of
nonnegative integers is said to be {\it unimodal\/} if there is an
index $m$ such that 
\beq
\label{a}
a_0\le a_1\le \ldots \le a_m\ge a_{m+1}\ge\ldots\ge a_r.
\eeq
Unimodal sequences arise in many aspects of combinatorics, geometry,
and algebra.  See the survey articles of Stanley~\cite{sta:lus} and
Brenti~\cite{bre:lus} for 
details.  Our interest is in the following well-known theorem.
\bth
\label{pkl}
Given $k,l\in\bbP$ the sequence
$$
p^{k,l}(0),p^{k,l}(1),\ldots,p^{k,l}(kl)
$$
is unimodal.\hqed
\eth
This result was first proved by Sylvester~\cite{syl:phu} using
invariant theory.  Since then, there have been a number of other
proofs.  In particular, Stanley~\cite{sta:wgh} derived this and much
more from the hard Lefschetz Theorem of algebraic geometry.
Proctor~\cite{pro:std} was able to reduce Stanley's proof to pure
linear algebra.  And finally, Kathy O'Hara~\cite{oha:ugc} gave a
combinatorial proof of this theorem.

\thicklines
\setlength{\unitlength}{2pt}
\bfi
\bpi(50,50)(-10,-10)
\put(-10,20){\makebox(0,0){$\la=$}}
\Gaaba \Gabdb \Gacdc \Gaded \Gaeee
\Gaaae \Gbabe \Gcbce \Gdbde \Gedee
\epi
\hspace{50pt}
\bpi(50,50)(-10,-10)
\put(-10,20){\makebox(0,0){$\ka=$}}
\Gaaba \Gabeb \Gacec \Gaddd \Gaede
\Gaaae \Gbabe \Gcbcc \Gcdce \Gdbdc \Gddde \Gebec
\epi
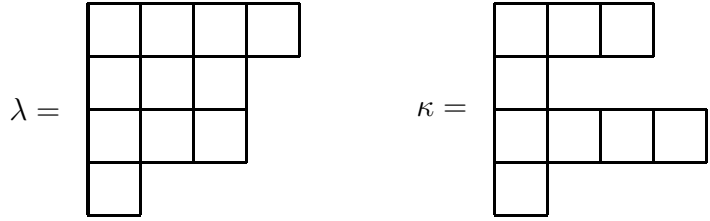
\capt{Ferrers diagrams for the partition $\la=(4,3,3,1)$ and the
  composition $\ka=(3,1,4,1)$}\label{laka}
\efi

A {\it composition $\ka$ of $n$\/}, written $\ka\cpn n$, is any
sequence $\ka=(\ka_1,\ldots,k_r)$ of positive integers summing to $n$.
Note that a composition need not be weakly decreasing.    All of the
definitions discussed so far have obvious analogues for compositions
so we will not bother restating them.  For example, the Ferrers
diagram of the composition $(3,1,4,1)$ is displayed in
Figure~\ref{laka}.  Although there is a large literature surrounding
partitions, composition have only recently aroused interest due to
their connection with quasi-symmetric functions~\cite{ehr:pha,ges:mpi},
the theory of patterns~\cite{hm:apl,sw:pac}, and the subword and
factor partial orders~\cite{bs:rmf,klrs:rif,sv:mfc}.

Let $c^{k,l}(n)$ be the number of compositions of $n$ fitting inside a
$k\times l$ rectangle.  In this note we will give a simple, algebraic proof
of the following conjecture of Vatter [personal communication].
\bth
Given $k,l\in\bbP$ the sequence
$$
c^{k,l}(0),c^{k,l}(1),\ldots,c^{k,l}(kl)
$$
is unimodal.
\eth

In the next section, we will prove this result by passing to the
generating function of the sequence.  The final section will include
some comments and an indication about how a combinatorial proof of
Vatter's conjecture might go.

\section{Unimodality of the composition sequence}

Let $a_0,a_1,\ldots,a_r$ be a sequence of real numbers and let $q$ be
a variable.  We consider the corresponding generating function
$f(q)=a_0+a_1q+\cdots a_r q^r$.  By convention, we let $a_i=0$ if
$i<0$ or $i>r$.  We will say that $f(q)$ has a given property if the
sequence itself does.

We will need the standard {\it $q$-analogue of $n$\/}, namely
$$
[n]=1+q+q^2+\cdots+q^{n-1}.
$$
It is well known that the generating function for the sequence
$p^{k,l}(n)$, $0\le n\le kl$, is the {\it $q$-binomial coefficient\/}
$$
\gau{k+l}{l}=\frac{[k+l]!}{[k]![l]!}
$$
where $[k]!=[k][k-1]\cdots[1]$.  So a restatement of Theorem~\ref{pkl}
is that the $q$-binomial coefficients are unimodal.

To prove the analogous result about compositions, we will need a
lemma.  It is not true, in general, that the product of two unimodal
polynomials is unimodal.  For example, if 
$f(q)=1+q+q^2+2.3q^3+2q^4$ then
$$
f(q)^2=1 + 2q + 3q^2 + 6.6q^3+ 9.6q^4 + 8.6q^5 + 9.29q^6 + 9.2q^7 + 4q^8.
$$
But we do have the following more specialized result.
\ble
\label{uni}
Let $f(q)$ be a unimodal polynomial and let $l\in\bbP$.  Then
$[l]f(q)$ is also unimodal. 
\ele
\prf
The lemma is clearly true if $l=1$, so assume $l\ge2$.
Suppose that $f(q)$ is the generating function for the
sequence~\ree{a}.  Also let $[l]f(q)=\sum_n b_n q^n$.  It follows
immediately from the definitions that
$$
b_0\le b_1\le\ldots \le b_m \qmq{and} 
b_{m+l-1}\ge b_{m+l}\ge\ldots\ge b_{r+l-1}.
$$
So the only way that $[l]f(q)$ could fail to be unimodal is if
$b_{i-1}>b_i<b_{i+1}$ for some $i$, $m<i<m+l-1$.  We will show that
the case $i=m+1$ leads to a contradiction as the other cases are
similar.

Suppose $b_m>b_{m+1}$ and $b_{m+1}<b_{m+2}$.  Expressing each $b_i$ in terms of the
$a_j$ and then cancelling terms gives $a_{m-l+1}>a_{m+1}$ and
$a_{m-l+2}<a_{m+2}$.  Using these inequalities and~\ree{a}, we have
$$
a_{m-l+1}>a_{m+1}\ge a_{m+2}>a_{m-l+2}.
$$
But this is a contradiction to~\ree{a} since $l\ge2$.\hqedm

Now let
$$
f^{k,l}=f^{k,l}(q)=\sum_{n\ge0} c^{k,l}(n) q^n.
$$
Our main result is as follows.
\bth
\label{fkl}
Let $k,l\in\bbP$.
\ben
\item[(a)] If $k\ge2$ then
$$
f^{k,l}=1+q[l]f^{k-1,l}.
$$
\item[(b)] The polynomial $f^{k,l}$ is unimodal.
\een
\eth
\prf
(a) Let $K^{k,l}$ be the set of compositions fitting inside a $k\times l$
rectangle, and let $K^{k,l}\spn{r}\sbe K^{k,l}$ be those compositions with first
part equal to $r$.  So we have the disjoint union 
\beq
\label{Kkl}
K^{k,l}=\{\ep\}\uplus\left(\biguplus_{r=1}^{l} K^{k,l}\spn{r}\right)
\eeq
where $\ep$ denotes the empty composition.  Removing the first part of
any $\ka\in K^{k,l}$ leaves a composition in $K^{k-1,l}$.  So
translating the union above into a generating function gives the
desired result.

(b)  We induct on $k$.  Clearly $f^{1,l}=[l+1]$, so we are done in the
base case.  If $k\ge2$ then using part (a) and the lemma finishes the
proof.\hqedm  

\section{Comments and open questions}

\subsection{Log concavity and symmetry}

A sequence $a_0,a_1,\ldots,a_r$ is {\it log concave\/} if
$a_i^2\ge a_{i-1} a_{i+1}$ for all $i$ with $0<i<r$.  The following
easily proved and well-known proposition gives a connection between
log-concavity and unimodality.
\bpr
Let $a_0,a_1,\ldots,a_r$ be a sequence of positive real numbers. 
If the sequence is log concave, then it is also unimodal.\hqed
\epr

Sometimes to prove a sequence is unimodal, it is actually easier to
prove that it satisfies the stronger log-concavity condition.  This is
because proving unimodality directly may involve finding the index
where the sequence is maximized, and that can be highly nontrivial.
However, the sequence $p^{k,l}(n)$, $0\le n\le kl$, is not
log concave in general.  So it should come as no surprise that neither
is $c^{k,l}(n)$, $0\le n\le kl$, and for much the same reason.  In
particular, if $k,l\ge2$ then both sequences start $1,1,2$ which
already violates log concavity.

Another common property of sequences is symmetry.  Say that
$a_0,a_1,\ldots,a_r$ is {\it symmetric\/} if $a_i=a_{r-i}$ for all
$i$, $0\le i\le r$.  By taking complements in the rectangle, it is
easy to see that $p^{k,l}(n)$, $0\le n\le kl$, is symmetric.  In
general, this property is not shared by compositions in a rectangle.
For example, if $k=l=2$ then the corresponding sequence is $1,1,2,2,1$.

\subsection{Lower order ideals}

Let $(P,\le)$ be a poset (partially ordered set).  Definitions for
terms from the theory of posets which are not given here can be found
in Stanley's book~\cite[Chapter 3]{sta:ec2}.  
A {\it lower order ideal\/} is $L\sbe P$ such that $x\in L$ and $y\le
x$ implies $y\in L$.  The {\it principal lower order ideal generated
by $x$\/} is the order ideal
$$
L(x)=\{y\in P\ |\ y\le x\}.
$$
Let $Y$ and $K$ denote Young's lattice and the poset of all
compositions, respectively.  Then the set of partitions in a rectangle
is the order ideal $Y(l^k)$ and similarly for compositions.

If $x,y\in P$ then {\it $x$ is covered by $y$\/}, written $x\lec y$,
if $x<y$ and there is no $z$ with $x<z<y$.
An {\it $x\hor y$ chain of length $n$\/} in $P$ is a subposet of the form
$x=x_0<x_1<\ldots<x_n=y$.  
This chain is {\it saturated\/} if each
inequality is actually a cover.
A poset is {\it graded\/} if it has a unique minimal element denoted
$\zh$, a unique maximal element denoted $\oh$, and every saturated
$\zh\hor\oh$ chain has the same length.  If $P$ is graded and $x\in P$
then all $\zh\hor x$ chains have the same length, called the 
{\it rank of $x$\/} and denoted $\rk x$.  In this case, the
{\it $n$th rank of $P$\/} is the subposet
$$
P_n=\{x\in P\ |\ \rk x = n\}.
$$

We will say that a graded poset $P$ has a property if the sequence of cardinalities
\beq
\label{P}
|P_0|, |P_1|,\ldots, |P_r|
\eeq
has that property, where $r=\rk \oh$.  We will sometimes preface the
property by ``rank-'' if clarification is needed.
So Theorem~\ref{pkl} can be restated as saying that the poset $Y(l^k)$
is unimodal.  It is natural to ask whether $Y(\la)$ is unimodal for
all partitions $\la$.  But this is too much to ask for, as
demonstrated by the following theorem of Stanton~\cite{sta:uyl}.
\bth[Stanton]
The lower order ideal $Y(8,8,4,4)$ is not unimodal.\hqed
\eth

In view of Stanton's result, it is perhaps surprising that all principal
lower order ideals in the composition poset $K$ are unimodal.
Given a graded poset $P$, we let $f^P=f^P(q)$ be the generating
polynomial for the sequence~\ree{P}.
The proof of the following theorem is so much like that of
Theorem~\ref{fkl} that we omit it.
\bth
Consider a composition $\ka\in K$.
\ben
\item[(a)]  Suppose $\ka=(\ka_1,\ldots,\ka_s)\neq\ep$, letting $l=\ka_1$ and
$\ga=(k_2,\ldots,k_s)$.  Then
$$
f^{K(\ka)}=1+q[l] f^{K(\ga)}.
$$
\item[(b)]  The polynomial $f^{K(\ka)}$ is unimodal.\hqed
\een
\eth

\subsection{A combinatorial proof?}

Theorem~\ref{fkl} is so easy to prove algebraically, one would think that
there is also an easy combinatorial proof.  But so far one has not
been found.  Here we present a possible inductive approach in
the hopes that someone else may be able to push it through.

Let $P$ be poset.  A {\it chain decomposition of $P$\/} is a family of
saturated chains $C_1,\ldots,C_a$  such that $P=\uplus_i C_i$.
If $P$ is graded  then we say an
$x\hor y$ chain in $P$ {\it symmetric\/} if $\rk y=\rk\oh-\rk x$.
A {\it symmetric chain decomposition\/} or {\it SCD\/} is a chain
decomposition where all the chains are symmetric.  It is easy to see
that if $P$ has an SCD then its rank sequence is symmetric and
unimodal.  

O'Hara~\cite{oha:ugc} constructed her ground-breaking
combinatorial proof of Theorem~\ref{pkl} as follows.
Let $Z(\la)$ be the poset of all partitions in $Y(\la)$ ordered by $\mu\le\nu$ if and
only if $|\mu|\le|\nu|$.  So for any partition $\la$, $Z(\la)$ has the
same set of ranks as does $Y(\la)$, but many more covering relations
in general.
\bth[O'Hara]
Given $k,l\in \bbP$, the poset $Z(l^k)$ has an SCD.\hqed
\eth
We note that it is still an open problem to give an SCD for $Y(l^k)$.

As mentioned above, $K(l^k)$ is not always
rank-symmetric.  But we can replace symmetry by another condition.
If $P$ is graded then we say that a chain decomposition is {\it
modal\/} (an {\it MCD\/}) if there is some rank $P_m$ such that every
$C_i$ contains an element of $P_m$.  
We call $P_m$ the {\it modal rank\/}.  The proof of the following
proposition is similar to the symmetric case, but we include it for
completeness.
\bpr
Let $P$ be a graded poset.  If $P$ has an MCD then $P$ is rank-unimodal.
\epr
\prf
Let $C_1,\ldots,C_a$ be an MCD and let $P_m$ be its modal rank.  We
will show that $|P_i|\le|P_{i+1}|$ for $i<m$ as the inequalities 
$|P_i|\ge|P_{i+1}|$ for $i\ge m$ are similar.  Let 
$P_i=\{x_1,\ldots,x_s\}$ and, since we have a cover, we can assume
that $C_j$ contains $x_j$ for $1\le j\le s$.  But each $C_j$ is
saturated and goes through rank $P_m$.  So for $1\le j\le s$, $C_j$
must contain an element $y_j$ in rank $P_{i+1}$.  By disjointness, the
$y_j$ are distinct and thus $|P_i|\le|P_{i+1}|$ as desired.\hqedm

We now ask the obvious questions.
\begin{question}
Does $K(l^k)$ have an MCD for all $k,l\in P$?  More generally, does
$K(\ka)$ have an MCD for all compositions $\ka$?
\end{question}
Note that the modal rank $P_m$ for $P=K(l^k)$ seems to occur when
$$
m=\ce{k(l+1)/2}-1
$$ where $\ce{\cdot}$ is the ceiling function.
Also note that there are other partial orders on the set of
compositions~\cite{bs:rmf,klrs:rif,sv:mfc} and they have the same set of
ranks as $K(l^k)$ (but not for a general $K(\ka)$).  Of these, the
partial order we are considering has the fewest covers.  So in may be
useful to consider one of the other orders instead.

It might be hoped that one could come up with an inductive description
of an MCD for $K(l^k)$ analogous to the inductive proof given of
Theorem~\ref{fkl}.  One possible way to do this is as follows.  For
simplicity, we will restrict ourselves to the case $l=2$.
Suppose we have an MCD $C_1,\ldots,C_a$ for $K(2^{k-1})$.  Then
using~\ree{Kkl} we can obtain a chain decomposition
$$
K(2^k)=\{\ep\}\uplus\{C_1',\ldots,C_a'\}\uplus\{C_1'',\ldots,C_a''\}
$$
where $C_i'$ is gotten by prefixing every element of $C_i$ by a one,
and $C_i''$ is obtained similarly using a two prefix.

Of course, this may be too many chains as not all of them will go
through rank $\ce{3k/2}-1$.  In particular, some of the $C_i'$ may be
too ``low'' and some of the $C_i''$ too ``high.''  (Also, $\ep$ must
be tacked onto some chain, but just use whichever $C_i'$ contains the
composition $(1)$.)  To rectify this, note that if $\ka'$ and $\ka''$
are the top elements of $C_i'$ and $C_i''$, respectively, then
by construction $\ka'\lec\ka''$.  So we can replace the pair
$C_i',C_i''$ by the pair $D_i',D_i''$ where
$$
D_i'=C_i'\uplus\{\ka''\}\qmq{and} D_i''=C_i''-\{\ka''\}.
$$
Note that this may result in $D_i''=\emp$ in which case we throw away
the chain.  Unfortunately, even with this correction the construction
breaks down when $k=9$ as some of the chains do not go through the
largest rank.  So some other modification will be needed to obtain an
MCD.

{\it Acknowledgement.}  I would like to thank Adam Goyt and Vince
Vatter for interesting discussions.

\bigskip
\bibliographystyle{acm}
\begin{small}
\bibliography{ref}
\end{small}

\end{document}